%% file: apprxmtly_M_dcmpsbl_sts.tex
\documentclass[11pt]{article}
\usepackage[dvipsnames]{xcolor}

\usepackage[utf8]{inputenc} 
\DeclareUnicodeCharacter{1EC7}{\d{\^{e}}}
\usepackage[T1]{fontenc} 

\usepackage[backend=biber, style=numeric, natbib=true, maxbibnames=10]{biblatex} 
\usepackage[autostyle=true]{csquotes} 

\usepackage{mathpazo} 
\usepackage{amssymb}
\usepackage{amsthm}
\usepackage{mathtools}
\usepackage{xifthen}
\usepackage[shortcuts]{extdash} 
\usepackage{hyperref}
\usepackage{enumitem}

\newtheorem{theorem}{Theorem}[section]
\newtheorem{proposition}[theorem]{Proposition}
\newtheorem{corollary}[theorem]{Corollary}
\newtheorem{lemma}[theorem]{Lemma}
\theoremstyle{definition}
\newtheorem{definition}[theorem]{Definition}
\newtheorem{remark}[theorem]{Remark}
\newtheorem{example}[theorem]{Example}

\NewDocumentCommand{\set}{mg}{\left\lbrace {#1} \IfValueT{#2}{\,\middle|\, {#2}} \right\rbrace}
\AtBeginDocument{\renewcommand{\matrix}[2][]{\ifthenelse{\isempty{#1}}{\begin{pmatrix}
        #2 \end{pmatrix}}{\left( \begin{array}{#1} #2 \end{array} \right)}}}

\newcommand{\R}{\mathbb{R}}
\newcommand{\N}{\mathbb{N}}
\newcommand{\T}{^{\mathsf{T}}}
\newcommand{\norm}[1]{\left\lVert {#1} \right\rVert}
\newcommand{\proj}[2][]{\ifthenelse{\isempty{#1}}{\pi\left({#2}\right)}{\pi_{#1}\left({#2}\right)}}
\newcommand{\dist}[2]{d\left({#1},{#2}\right)}
\newcommand{\hd}[2]{d_{\mathsf{H}}\left({#1},{#2}\right)} 
\newcommand{\rc}[1]{0^{+}{#1}}
\newcommand{\polar}[1]{{#1}^{\circ}}

\def\fchar#1#2@{#1} 
\def\estr#1#2@{#2} 
\NewDocumentCommand{\ccc}{m}{{\fchar #1@}^*{\estr #1@}}
\NewDocumentCommand{\cc}{m}{\expandafter\ccc\expanded{{#1}}} 
\NewDocumentCommand{\biccc}{m}{{\fchar #1@}^{**}{\estr #1@}}
\NewDocumentCommand{\bicc}{m}{\expandafter\biccc\expanded{{#1}}} 

\renewcommand{\geq}{\geqslant}
\renewcommand{\leq}{\leqslant}

\DeclareMathOperator{\cl}{cl}
\let\int\relax
\DeclareMathOperator{\int}{int}
\DeclareMathOperator{\relint}{relint}
\DeclareMathOperator{\dom}{dom}

\DeclareMathOperator{\cone}{cone}

\DeclareUnicodeCharacter{1EE5}{\d{u}}

\newcommand{\figpath}{figures/}

\begin{document}

\title{Convex sets approximable as the sum of a compact set and a
  cone}

\author{Daniel Dörfler\thanks{daniel.doerfler@uni-jena.de, Friedrich Schiller University Jena, Germany}
  \and Andreas Löhne\thanks{andreas.loehne@uni-jena.de, Friedrich Schiller University Jena, Germany}}

\date{\today}

\maketitle

\begin{center}
  \textit{Dedicated to the memory of Professor Alfred Göpfert.}
  \vspace{.5cm}
\end{center}

\begin{abstract}
  The class of convex sets that admit approximations as Minkowski sum
  of a compact convex set and a closed convex cone in the Hausdorff
  distance is introduced. These sets are called approximately
  Motzkin\-/decomposable and generalize the notion of
  Motzkin\-/decomposability, i.e. the representation of a set as the sum
  of a compact convex set and a closed convex cone. We characterize
  these sets in terms of their support functions and show that they
  coincide with hyperbolic sets, i.e. convex sets contained in the sum
  of their recession cone and a compact convex set, if their recession
  cones are polyhedral but are more restrictive in general. In
  particular we prove that a set is approximately Motzkin-decomposable
  if and only if its support function has a closed domain relative to
  which it is continuous.
\end{abstract}

\input{\secpath sec1}

\input{\secpath sec2}

\input{\secpath sec3}

\input{\secpath sec4}

\vskip 6mm
\noindent\textbf{Acknowledgements}
\noindent The authors thank Juan Enrique Martínez-Legaz for useful
comments on the contents of the article, in particular, for pointing
them to \cite{Bai83} and the notion of hyperbolic convex set.

\sloppy
\printbibliography[heading=bibintoc]
\fussy

\end{document}

%% file: sections/sec1.tex
\section{Introduction}\label{sec_1}

\begin{figure}[h!]
  \centering
  \includegraphics[]{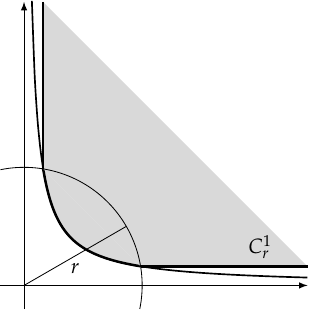}
  \hspace{.5cm}
  \includegraphics[]{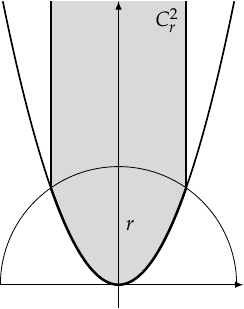}
  \caption{When can $C$ be approximated by truncations?}
  \label{fig_1}
\end{figure}

Consider the closed convex sets $C^1 = \set{x \in \R^2}{x_1x_2 \geq 1,
x_1 \geq 0}$ and $C^2 = \set{x \in \R^2}{x_2 \geq x_1^2}$ with their
recession cones $\rc{C^1} = \set{x \in \R^2}{x \geq 0}$ and $\rc{C^2}
= \set{x \in \R^2}{x_1 = 0, x_2 \geq 0}$. For $r > 0$ we understand a
truncation $C_r^i$ of these sets to be $C_r^i = \left(C^i \cap
rB\right) + \rc{C^i}$, $i=1,2$, where $B$ denotes the Euclidean unit
ball. That is, we restrict $C^i$ to a compact subset of itself given
by a ball of radius $r$ and add the recession cone afterwards,
yielding a subset of $C^i$. This situation is depicted in Figure
\ref{fig_1}. Clearly, the truncations $C_r^1$ of $C^1$ converge to
$C^1$ with respect to the Hausdorff distance if $r$ tends to
$+\infty$. However, the same is not true for $C^2$ because, for every
$r > 0$, the Hausdorff distance between $C_r^2$ and $C^2$ is 
infinite. Note that there is a geometric difference between $C^1$ and
$C^2$: $C^1 \subseteq \set{0} + \rc{C^1}$, but there does not exist a
compact set $M$ satisfying $C^2 \subseteq M + \rc{C^2}$. In view of
this property, $C^1$ is called hyperbolic, cf. \cite{Bai83,
  MLRR93}. In this work we investigate the following questions
regarding a closed convex set $C \subseteq \R^n$:

\begin{enumerate}
\item What characterizes the property that $C$ can be approximated by
  truncations $C_r$ in the above sense?
\item How is this property related to hyperbolicity of $C$? Do both
  properties coincide (under certain assumptions)?
\end{enumerate}

This is motivated by the observation that conducting analysis on $C$
presents various challenges in the presence of unboundedness. The
above example demonstrates that the Hausdorff distance is only
limitedly suitable for the comparison of unbounded sets and the
recession cone may only partly capture the asymptotic behavior of
$C$. The latter is evident from the fact that $\rc{C} = \rc{C_r}$
whenever $C_r$ is nonempty. Hence, from a practical perspective, it
seems favourable to reduce investigations to a bounded subset of $C$.

What needs to be adressed is whether it is possible to select a subset
a priori in such a way as to not neglect important information about
$C$. One class of convex sets that permit such a choice are called
M-decomposable sets and have been studied by various authors over the
last decade \cites{GGMT10, GMT10, GIMT13, Sol20}. These are the sets
that can be represented as the Minkowski sum of a compact convex set
and a closed convex cone. They have been characterized in different
ways, e.g. in terms of support functions and via an associated vector
optimization problem and its set of efficient points \cite{GGMT10},
with regard to their extreme points \cite{GMT10} and in the sense of
intersections of sets with certain hyperplanes
\cite{GIMT13}. Nevertheless, imposing M-decomposability on $C$ is a
strong assumption. For example, both of the sets $C^1$ and $C^2$ above
are not M-decomposable. However, their truncations always are. Hence,
the first question might be restated as:

\begin{enumerate}
\item[$1^*$.] When can a closed convex set $C$ be approximated by
  M-decomposable sets?
\end{enumerate}

To address these questions we introduce and study the class of
approximately M-decomposable sets as those closed convex sets that can
be approximated by M-decomposable sets in the Hausdorff distance. Our
main result is a characterization of this class in terms of support
functions and their demarcation from M-decomposable and hyperbolic
sets. In the following section we provide preliminary definitions and
notation. Section \ref{sec_3} contains the definition of approximately
M-decomposable sets and the main results of this manuscript. In
particular we show that approximately M-decomposable sets fit between
M-decomposable and hyperbolic sets. They coincide with the latter
given their recession cones are polyhedral. Moreover, we prove that a
set is approximately M-decomposable if and only if its support
function has a closed domain relative to which it is continuous.

%% file: sections/sec2.tex
\section{Preliminaries}\label{sec_2}

Given sets $M$, $P \subseteq \R^n$ and a scalar $\alpha \in \R$ we
denote by $\cl M$, $\int M$, $\relint M$, $M+P$ and $\alpha M$ the
\emph{closure}, \emph{interior}, \emph{relative interior} of $M$, the
\emph{Minkowski sum} of $M$ and $P$ and the \emph{dilation} of $M$ by
$\alpha$, respectively. The \emph{Euclidean norm} of a vector $x \in
\R^n$ and the \emph{Euclidean unit ball} are expressed as $\norm{x}$
and $B$, respectively. A set $C \subseteq \R^n$ is called
\emph{convex} if for every $x, y \in C$ and scalar $\lambda \in
[0,1]$, $\lambda x + (1-\lambda) y \in C$ holds. A \emph{cone} is a
set $C$ such that for every $x \in C$ and $\mu \geq 0$ it holds $\mu x
\in C$. The \emph{recession cone} of a convex set $C \subseteq \R^n$
is defined as the set $\set{d \in \R^n}{\forall\; x \in C, \mu \geq 0
\colon x + \mu d \in C}$ and denoted by $\rc{C}$. If $C$ is closed, so
is $\rc{C}$. The \emph{polar} $\polar{C}$ of $C$ is the set of linear
functions that are bounded above on $C$ by $1$, i.e. $\polar{C} =
\set{y \in \R^n}{\forall\; x \in C \colon y\T x \leq 1}$. If $C$ is a
cone, then $\polar{C} = \set{y \in \R^n}{\forall\; x \in C \colon y\T
x \leq 0}$. Note that $\polar{C}$ is always closed. Given nonempty
sets $C^1$, $C^2 \in \R^n$ the \emph{Hausdorff distance}
$\hd{C^1}{C^2}$ between $C^1$ and $C^2$ is defined as
\begin{equation*}
  \hd{C^1}{C^2} = \max\set{\adjustlimits\sup_{x \in C^1}\inf_{y \in
      C^2}\norm{x-y},\; \adjustlimits\sup_{x \in C^2}\inf_{y \in
      C^1}\norm{x-y}}.
\end{equation*}
Note that $\hd{C^1}{C^2}$ may be infinite if any of the sets is
unbounded.  However, it is well known that $\hd{\cdot}{\cdot}$ defines
a metric on the class of nonempty compact subsets of $\R^n$
\cite{Hau49}. Moreover, it can equivalently be expressed as
\begin{equation}\label{eq_2}
  \hd{C^1}{C^2} = \inf\set{\varepsilon > 0}{C^1 \subseteq C^2 +
    \varepsilon B,\; C^2 \subseteq C^1 + \varepsilon B},
\end{equation}
see \cite{Gru07}.

A function $f \colon \R^n \to \R \cup \set{+\infty}$ is called
\emph{convex} if $f(\lambda x + (1-\lambda) y) \leq \lambda f(x) +
(1-\lambda) f(y)$ holds for every $x$, $y \in \R^n$ and $\lambda \in
[0,1]$. The \emph{domain} $\dom f$ of $f$ is the set of points at
which $f$ is finite, i.e. $\dom f = \set{x \in \R^n}{f(x) <
+\infty}$. To every $f$ we assign a \emph{conjugate} function $\cc{f}$
defined by $\cc{f}(y) = \sup\set{y\T x - f(x)}{x \in \R^n}$. The
\emph{support function} $\sigma_C$ of a nonempty closed convex set $C$
is given as $\sigma_C(d) = \sup\set{d\T x}{x \in C}$. Its conjugate is
called \emph{indicator function} of $C$ and denoted $\delta_C$.

%% file: sections/sec3.tex
\section{Approximately M-decomposable sets}\label{sec_3}

In this section we give an answer to the questions posed in the
introduction. We understand a truncation of $C$ to be the sum of a
compact subset of $C$ and $\rc{C}$ in the following sense.

\begin{definition}
  Let $C \subseteq \R^n$ be nonempty closed and convex. The set $C_r =
  \left(C \cap rB\right) + \rc{C}$ is called a \emph{truncation} of
  $C$ of radius $r > 0$.
\end{definition}

Clearly, $C_r \subseteq C$ for all $r > 0$. Studying sets that are
approximately trunctations of themselves is motivated by the fact
that, given such a set, it suffices to conduct analyses, such as
minimizing some function over the set, on a compact subset to obtain
good results in an approximate sense. To this end we introduce the
class of approximately M(otzkin)-decomposable convex sets.

\begin{definition}\label{def_3.1}
  A nonempty closed convex set $C \subseteq \R^n$ is called
  \emph{approximately M-decomposable} if, for every $\varepsilon > 0$,
  there exists $r > 0$ such that $C \subseteq C_r + \varepsilon B$.
\end{definition}

By Equation \eqref{eq_2}, Definition \ref{def_3.1} is equivalent to
the existence of a truncation $C_{r(\varepsilon)}$ such that the
Hausdorff distance between $C_{r(\varepsilon)}$ and $C$ is at most
$\varepsilon$. The designation approximately M-decomposable is derived
from the stronger property of M-decomposability.

\begin{definition}[{see \cite{GGMT10}}]\label{def_3.2}
    A nonempty closed convex set $C \subseteq \R^n$ is called
    \emph{M(otzkin)-decomposable} if there exists a compact convex set
    $M \subseteq \R^n$ such that $C = M + \rc{C}$. Such set $M$ is
    called a compact component of $C$.
\end{definition}

If $C$ is closed and convex, then every nonempty trunctation $C_r$ of
$C$ is M\-/decomposable with compact component $C \cap rB$. Moreover,
if $C$ is M-decomposable, there exist truncations of $C$ that coincide
with $C$. If $M$ is any compact component of $C$, then $C_r = C$ for
all $r \geq \max\set{\norm{x}}{x \in M}$. In particular, $C$ is
approximately M-decomposable in this case. A relaxation of the
condition of M-decomposability leads to hyperbolic sets.

\begin{definition}[{cf. \cites{Bai83, MLRR93}}]\label{def_3.3}
  A nonempty closed convex set $C \subseteq \R^n$ is called
  \emph{hyperbolic} if there exists a compact convex set $M
  \subseteq \R^n$ such that $C \subseteq M + \rc{C}$.
\end{definition}

We point out that hyperbolic sets also exist under the terms
$\rc{C}$-bounded \cite{Luc89} and self-bounded \cite{Doe22} sets in
the literature. Moreover, it is known that $C$ is hyperbolic if and
only if $\dom\sigma_C$ is closed, see \cite[Proposition 5]{Bai83}.

Clearly, every M-decomposable set $C$ is also hyperbolic, where the
set $M$ in Definition \ref{def_3.3} may be chosen as any compact
component of $C$, but the converse does not hold as is seen from the
example $C^1 = \set{x \in \R^2}{x_1x_2 \geq 1, x_1 \geq 0}$ and the
fact $\rc{C^1} = \set{x \in \R^2}{x \geq 0}$, cf. Section \ref{sec_1}
and Figure \ref{fig_1}. The set $C^1$ is the epigraph of the function
$x \mapsto x^{-1}$ restricted to the nonnegative real line. It is not
difficult to see that this set is approximately M-decomposable. For
$\varepsilon > 0$ set $r^2 \geq \varepsilon^2 +
\varepsilon^{-2}$. Example \ref{ex_3.12} below shows that not every
hyperbolic set is approximately M-decomposable. However, the converse
statement holds, i.e. the class of hyperbolic sets is larger in
general.

\begin{proposition}\label{prop_3.5}
  Let $C \subseteq \R^n$ be approximately M-decomposable. Then $C$ is
  hyperbolic.
\end{proposition}
\begin{proof}
  Let $\varepsilon > 0$. Then there exists $r > 0$ such that $C
  \subseteq C_r + \varepsilon B$. Now,
  \[
    C_r + \varepsilon B = \left(C \cap rB\right) + \varepsilon B +
    \rc{C}
  \]
  and $\left(C \cap rB\right) + \varepsilon B$ is compact and
  convex. Thus, $C$ is hyperbolic.
\end{proof}

In order to characterize approximately M-decomposable sets we need
another expression for the Hausdorff distance between certain closed
convex sets. The following is a generalization of a well-known result
for compact convex sets, see e.g. \cite[Proposition 6.3]{Gru07}.

\begin{proposition}\label{prop_3.6}
  Let $C^1$, $C^2 \subseteq \R^n$ be nonempty closed convex sets whose
  support functions have identical domains $D$, i.e. $D =
  \dom\sigma_{C^1} = \dom\sigma_{C^2}$. Then
  \[
    \hd{C^1}{C^2} = \sup_{d \in D \cap B}\left|\sigma_{C^1}(d) -
      \sigma_{C^2}(d)\right|.
  \]
  Furthermore, $\hd{C^1}{C^2} = +\infty$ whenever $\dom\sigma_{C^1}
  \neq \dom\sigma_{C^2}$.
\end{proposition}
\begin{proof}
  Since $C^1$ and $C^2$ are nonempty, $D$ is also nonempty because
  $\sigma_{C^1}(0) = \sigma_{C^2}(0) = 0$, i.e. $0 \in D$. Denote by
  $\dist{x}{C^1}$ the Euclidean distance from $x$ to $C^1$,
  i.e. $\dist{x}{C^1} = \inf_{y \in C^1} \norm{x-y}$. Then
  $\hd{C^1}{C^2}$ can be written as $\max\set{\sup_{x \in
      C^1}\dist{x}{C^2},\; \sup_{x \in C^2}\dist{x}{C^1}}$. Moreover, it
  holds $\dist{x}{C^1} = \inf_{y \in \R^n}\left(\norm{x-y} +
    \delta_{C^1}(y)\right)$. According to \cite[Theorem 16.4]{Roc70} the
  conjugate of $\dist{\cdot}{C^1}$ is thus given as
  $\cc{\dist{\cdot}{C^1}} = \delta_B + \sigma_{C^1}$. Since
  $\dist{\cdot}{C^1}$ is closed proper and convex it holds
  \[
    \begin{aligned}
      \dist{x}{C^1} &= \bicc{\dist{x}{C^1}}\\
      &= \sup_{d \in \R^n}\left(x\T d - \delta_B(d) -
        \sigma_{C^1}(d)\right)\\
      &= \sup_{d \in D \cap B}\left(x\T d - \sigma_{C^1}(d)\right).
    \end{aligned}
  \]
  The first equation is due to \cite[Theorem 12.2]{Roc70} and the
  third holds because $\dom\sigma_{C^1} = D$. Analogously one has
  $\dist{x}{C^2} = \sup_{d \in D \cap B}\left(x\T d -
    \sigma_{C^2}(d)\right)$. Combining both expressions yields
  \[
    \begin{aligned}
      \hd{C^1}{C^2} &= \max\set{\adjustlimits\sup_{x \in C^2}\sup_{d
          \in D \cap B}\left(x\T d - \sigma_{C^1}(d)\right),\;
        \adjustlimits\sup_{x \in C^1}\sup_{d \in D \cap B}\left(x\T d
          - \sigma_{C^2}(d)\right)}\\
      &= \sup_{d \in D \cap B}\max\set{\sigma_{C^2}(d) -
        \sigma_{C^1}(d),\; \sigma_{C^1}(d) - \sigma_{C^2}(d)}\\
      &= \sup_{d \in D \cap B}\left|\sigma_{C^1}(d) -
        \sigma_{C^2}(d)\right|.
    \end{aligned}
  \]
  For the second part assume, without loss of generality, $d \in
  \dom\sigma_{C^1}\setminus\dom\sigma_{C^2}$, $\norm{d} = 1$. Then for
  every $k \in \N$ there exists $y_k \in C^2$, such that $d\T y_k >
  k$. Denote by $H$ the set $\set{x \in \R^n}{d\T x \leq
    \sigma_{C^1}(d)}$. Since $C^1 \subseteq H$, it holds
  $\dist{y_k}{C^1} \geq \dist{y_k}{H} = d\T y_k - \sigma_{C^1}(d) > k
  - \sigma_{C^1}(d)$ for all $k \geq \sigma_{C^1}(d)$. Taking the
  limit $k \to +\infty$ yields the result.
\end{proof}

\begin{remark}
  Note that, even if $\dom\sigma_{C^1} = \dom\sigma_{C^2}$, it is
  possible, that $\hd{C^1}{C^2} = +\infty$. Consider, for example, the
  sets $C^1 = \set{x \in \R^2}{x_2 \geq x_1^2}$ and $C^2 = 2C^1$. The
  domain of their support functions is the non-closed set $\set{x \in
    \R^2}{x_2 < 0} \cup \set{0}$. For $d_k = \left(1,\; -1/k\right)\T$ one
  has $\sigma_{C^1}(d_k) = k/2$ and $\sigma_{C^2}(d_k) = k$, i.e.
  \[
    \sigma_{C^2}\left(\frac{d_k}{\norm{d_k}}\right) -
    \sigma_{C^1}\left(\frac{d_k}{\norm{d_k}}\right) =
    \frac{k^2}{2\sqrt{k^2+1}}
  \]
  which is unbounded in $k \in \N$.
\end{remark}

We have already seen that not every closed convex set $C$ is
approximately M-decomposable, i.e. can be approximated by truncations
in the Hausdorff distance. However, truncations always approximate $C$
in a weaker sense.

\begin{definition}[{cf. \cite{RW98}}]
  A sequence of $\set{C^k}_{k \in \N}$ of closed subsets of $\R^n$ is
  said to \emph{converge in the sense of Painlevé-Kuratowski} or
  \emph{PK-converge} to a closed set $C \subseteq \R^n$ if
  \[
    C = \set{x \in \R^n}{\lim_{k \to +\infty}\dist{x}{C^k} = 0}.
  \]
\end{definition}

For bounded sequences PK-convergence and convergence with respect to
the Hausdorff distance coincide, but the concepts are distinct in
general, see \cite{RW98}.

\begin{lemma}\label{lem_3.9}
  Let $C \subseteq \R^n$ be nonempty closed convex. Then $\{C_r\}_{r
    \in \N}$ PK-converges to $C$.
\end{lemma}
\begin{proof}
  Let $x \in C$. Then $x \in C_r$ for all $r \geq \norm{x}$. Hence,
  $\lim_{r \to +\infty}\dist{x}{C_r}$ exists and is zero. On the
  contrary, let $x \notin C$. Since $C$ is closed, $\dist{x}{C} >
  0$. Moreover, $\dist{x}{C_r} \geq \dist{x}{C}$ holds for all $r \in
  \N$ because $C_r \subseteq C$. Thus $x \notin \set{x \in \R^n}{\lim_{k
      \to +\infty}\dist{x}{C^k} = 0}$.
\end{proof}

We are now ready to prove the main result of this work.

\begin{theorem}\label{thm_3.10}
  Let $C \subseteq \R^n$ be a nonempty closed convex set. The
  following are equivalent:
  \begin{enumerate}[label=(\roman*)]
  \item $C$ is approximately M-decomposable,
  \item\label{thm_3.10ii} $\dom\sigma_C$ is closed and $\sigma_C$ is
    continuous relative to $\dom\sigma_C$.
  \end{enumerate}
  Moreover, the set $\dom\sigma_C$ in \ref{thm_3.10ii} equals
  $\polar{\left(\rc{C}\right)}$.
\end{theorem}
\begin{proof}
  Assume $C$ is approximately M-decomposable, i.e.
  \begin{equation}\label{eq_3}
    \forall\;\varepsilon > 0\;\exists\;r \geq 0\colon C_r \subseteq C
    \subseteq C_r + \varepsilon B.
  \end{equation}
  For the remainder of the proof assume $\varepsilon$ and $r$ are
  fixed such that the inclusions in \eqref{eq_3} hold. In particular,
  $C_r + \varepsilon B$ and, consequently, $C_r$ are nonempty. It
  follows
  \begin{equation}\label{eq_4}
    \sigma_{C_r} \leq \sigma_C \leq \sigma_{C_r + \varepsilon B}
  \end{equation}
  for the support functions of $C_r$, $C$ and $C_r + \varepsilon B$,
  see \cite[Corollary 13.1.1]{Roc70}. Since all three sets have the
  same recession cone $\rc{C}$, it follows from \cite[Corollary
  14.2.1]{Roc70} that
  \[
    \cl\dom\sigma_{C_r} = \cl\dom\sigma_C = \cl\dom\sigma_{C_r +
      \varepsilon B} = \polar{\left(\rc{C}\right)}.
  \]
  However, for $d \in \polar{\left(\rc{C}\right)}$ one has
  \[
    \begin{aligned}
      \sigma_{C_r + \varepsilon B}(d) = \sigma_{\left(C \cap rB\right)
        + \varepsilon B}(d) < +\infty.
    \end{aligned}
  \]
  Thus $\dom\sigma_{C_r} = \dom\sigma_C = \dom\sigma_{C_r +
    \varepsilon B} = \polar{\left(\rc{C}\right)}$ due to Inequality
  \eqref{eq_4}.

  The support functions $\sigma_{C_r}$ and $\sigma_{C \cap rB}$
  coincide on $\polar{\left(\rc{C}\right)}$. Since $C \cap rB$ is
  compact, its support function is finite everywhere and therefore
  continuous according to \cite[Corollary 10.1.1]{Roc70}. Thus,
  $\sigma_{C_r}$ is continuous relative to its domain
  $\polar{\left(\rc{C}\right)}$. Using the expression for the
  Hausdorff distance from Equation \eqref{eq_2}, \eqref{eq_3} is
  equivalent to
  \begin{equation}\label{eq_5}
    \forall\;\varepsilon > 0\;\exists\;r \geq 0\colon \hd{C_r}{C} \leq
    \varepsilon.
  \end{equation}
  Since
  \begin{equation}\label{eq_6}
      r \leq \bar{r} \implies C_r \subseteq C_{\bar{r}} \subseteq C,
  \end{equation}
  Statement \eqref{eq_5} is equivalent to $\lim_{r \to
    +\infty}\hd{C_r}{C} = 0$. Applying Proposition \ref{prop_3.6}
  gives
  \begin{equation}\label{eq_7}
    \adjustlimits\lim_{r \to +\infty}\sup_{d \in
      \polar{\left(\rc{C}\right)} \cap B}\left|\sigma_C(d) -
      \sigma_{C_r}(d)\right| = 0.
  \end{equation}
  Thus, $\sigma_{C_r}$ converges uniformly to $\sigma_C$ on
  $\polar{\left(\rc{C}\right)}$. As $\sigma_{C_r}$ is continunous on its
  domain, the Uniform Limit Theorem \cite[Theorem 21.6]{Mun00} yields
  that $\sigma_C$ is continuous on its domain as well.

  Now, assume $\dom\sigma_C$ is closed and $\sigma_C \colon
  \dom\sigma_C \to \R$ is continuous. Again, $\dom\sigma_C =
  \polar{\left(\rc{C}\right)}$ by \cite[Corollary 14.2.1]{Roc70}. We
  show that $\sigma_{C_r}$ converges pointwise to $\sigma_{C}$ on
  $\polar{\left(\rc{C}\right)}$ as $r \to +\infty$. Property
  \eqref{eq_6} implies that $\lim_{r \to +\infty}\sigma_{C_r}(d)$ exists
  for all $d \in \polar{\left(\rc{C}\right)}$ and is finite. In
  particular it holds $\lim_{r \to +\infty}\sigma_{C_r}(d) \leq
  \sigma_C(d)$. Assume $\sigma_{C_r}$ does not converge to pointwise to
  $\sigma_C$ on $\polar{\left(\rc{C}\right)}$, i.e. there exists
  $\bar{d} \in \polar{\left(\rc{C}\right)}$ for which the inequality
  holds strictly. This implies $\bar{d} \neq 0$. Let $\gamma = \lim_{r
    \to +\infty}\sigma_{C_r}(\bar{d})$. Then $C_r \subseteq H := \set{x
    \in \R^n}{\bar{d}\T x \leq \gamma}$ for every $r \in \N$ and there
  exists $\bar{x} \in C$ such that $\bar{d}\T \bar{x} > \gamma$. It
  holds
  \[
    \dist{\bar{x}}{C_r} \geq \dist{\bar{x}}{H} = \frac{\bar{d}\T
      \bar{x} - \gamma}{\norm{\bar{d}}} > 0.
  \]
  This is a contradiction because $\set{C_r}_{r \in \N}$ PK-converges
  to $C$ by Lemma \ref{lem_3.9}, i.e. $C = \set{x \in \R^n}{\lim_{r \to
      +\infty}\dist{x}{C_r} = 0}$. Thus $\sigma_{C_r}$ converges
  pointwise to $\sigma_C$ on $\polar{\left(\rc{C}\right)}$.

  Since $\sigma_C$ is continuous and $\sigma_{C_r}$ is monotonically
  non-decreasing in $r$, Dini's theorem \cite[Theorem 12.1]{Jos05}
  yields that the convergence is uniform on $\polar{\left(\rc{C}\right)}
  \cap B$, i.e. Equation \eqref{eq_7} holds.  Finally, using the
  equivalence between \eqref{eq_7}, \eqref{eq_5} and \eqref{eq_3}
  proves that $C$ is approximately M-decomposable.
\end{proof}

\begin{remark}
  In the above proof we have shown the pointwise convergence of the
  support functions of truncations $C_r$ to the support function of
  $C$ on $\polar{\left(\rc{C}\right)}$ using the PK-convergence of
  $\set{C_r}_{r \in \N}$. We point out that this may also be achieved
  using the different concept of $\mathcal{C}$-convergence introduced
  in \cite{LZ06}. To this end one needs to take into account that
  $\sigma_C = \sup_{r \in \N}\sigma_{C_r}$ and apply \cite[Theorem
  5.9]{LZ06} and \cite[Proposition 5.3 (iv)]{LZ06}.
\end{remark}

\begin{remark}
  In \cite{GK59} \textit{continuous convex sets} are introduced, which
  are similar to, but different from, approximate M-decomposable sets. A
  convex set is called continuous if its support function is continuous
  everywhere (and not necessarily only relative to its domain). For
  example, the set $C^1$ from Section \ref{sec_1} is approximately
  M-decomposable, but not continuous, whereas $C^2$ is continuous, but
  the domain of its support function is not closed. Thus it is not
  approximately M-decomposable, see \cite[p. 1]{GK59}. Clearly, every
  compact convex set is both continuous and (approximately)
  M-decomposable.
\end{remark}

Using Theorem \ref{thm_3.10} we can show that hyperbolicity does not
imply approximate M-decomposability.

\begin{example}\label{ex_3.12}
  Let $C^2 = \set{x \in \R^2}{x_2 \geq x_1^2}$, cf. Section
  \ref{sec_1} and Figure \ref{fig_1}, and consider the following set
  \[
    C = \left(C^2 \times \set{0}\right) + \cl\cone\left(C^2 \times
      \set{1}\right)
  \]
  in $\R^3$. By \cite[Theorem 8.2]{Roc70} it holds
  \[
    \cl\cone\left(C^2 \times \set{1}\right) = \cone\left(C^2 \times
      \set{1}\right) \cup \left(\rc{C^2} \times \set{0}\right).
  \]
  Now, according to \cite[Corollary 9.1.2]{Roc70} $C$ is closed and
  one has $\rc{C} = \cl\cone\left(C^2 \times \set{1}\right)$. To see
  that $C$ is hyperbolic, choose any $x \in C$. Then $x = c + s$ for
  some $c \in C^2 \times \set{0}$ and $s \in \rc{C}$. For $d =
  (0,0,-1)\T$ it holds $x = d+(c-d)+s$. Since $c-d \in \rc{C}$, $c-d+s
  \in \rc{C}$ as well. Hence, $C \subseteq \set{d} + \rc{C}$.

  According to \cite[Theorem 3.1]{BBW23} one has
  \[
    \polar{\left(\rc{C}\right)} = \cl\cone\left(\polar{C^2} \times
      \set{-1}\right)
  \]
  and a simple calculation shows that $\polar{C^2} = \set{x \in
    \R^2}{x_2 \leq -x_1^2/4}$. Therefore,
  \[
    d_n = \matrix{n^{-1}\\-(2n)^{-2}\\-1} \in
    \polar{\left(\rc{C}\right)}
  \]
  for every $n \in \N$ and $\lim_{n \to +\infty}d_n = d$. Now,
  \[
    \sigma_C(d_n) \geq d_n\T \matrix{n\\n^2\\0} = \frac{3}{4} > 0.
  \]
  However, $\sigma_C(d) \leq 0$ because for $x \in C$ it holds $x_3
  \geq 0$. Hence, $\sigma_C$ is not continuous at $d \in
  \polar{\left(\rc{C}\right)}$ and by Theorem \ref{thm_3.10} not
  approximately M-decomposable.
\end{example}

If an additional assumption is made, hyperbolicity and approximate
M-decomposability are equivalent.

\begin{corollary}
  Let $C \subseteq \R^n$ be nonempty closed convex and $\rc{C}$ be
  polyhedral. Then $C$ is approximately M-decomposable if and only if
  it is hyperbolic.
\end{corollary}
\begin{proof}
  By Proposition \ref{prop_3.5}, the first implication is true
  regardless of the polyhedrality of $\rc{C}$.

  Assume that $C$ is hyperbolic, i.e. $C \subseteq M + \rc{C}$ for
  some compact convex set $M$. This implies $\dom\sigma_C \subseteq
  \polar{\left(\rc{C}\right)}$ and thus actually $\dom\sigma_C =
  \polar{\left(\rc{C}\right)}$ by \cite[Corollary 14.2.1]{Roc70}. Since
  $\rc{C}$ is polyhedral, so is $\polar{\left(\rc{C}\right)}$. Hence it
  is locally simplicial, see \cite[Theorem 20.5]{Roc70}, and
  \cite[Theorem 10.2]{Roc70} implies that $\sigma_C$ is continuous
  relative to its domain. Now, apply Theorem \ref{thm_3.10}.
\end{proof}

The corollary implies that, whenever $C \subseteq \R^2$, it is
approximately M-decomposable if and only if it is hyperbolic because
every closed convex cone in $\R^2$ is polyhedral. Thus, the set from
Example \ref{ex_3.12} is in the lowest possible dimension to
demonstrate that hyperbolicity is, in general, not equivalent to
approximate M-decomposability.

%% file: sections/sec4.tex
\section{Conclusion}

In this article we have introduced the class of approximately
M-decomposable sets as those closed convex sets that can be
approximated arbitrarily well by truncations with respect to the
Hausdorff distance. We have shown that they generalize the class of
M-decomposable sets and are a special case of hyperbolic sets but
coincide with neither in general as demonstrated by
examples. Furthermore, we characterized the approximately
M-decomposable sets in terms of their support functions. Finally we
proved that, when considering only polyhedral recession cones,
approximate M-decomposability and hyperbolicity are equivalent.